\documentclass[iccmp]{ipbook}

\RequirePackage{amsthm,amsmath,amssymb}
\RequirePackage{hyperref}
\hypersetup{
    colorlinks,
    linktocpage,
    linkcolor={blue}
    }
\RequirePackage{tabularx}

\usepackage[all]{xy}    
\usepackage{amsmath,amsthm,amssymb,mathrsfs,mathtools,bm,eucal,tensor}                  

\usepackage{bbm}

\usepackage{fancyhdr}
\newcommand{\BA}{{\mathbb {A}}}

\newcommand{\BC}{{\mathbb {C}}}

\newcommand{\BG}{{\mathbb {G}}}
\newcommand{\BH}{{\mathbb {H}}}

\newcommand{\BR}{{\mathbb {R}}}

\newcommand{\BX}{{\mathbb {X}}}
\newcommand{\BY}{{\mathbb {Y}}}
\newcommand{\BZ}{{\mathbb {Z}}}

\newcommand{\CC}{{\mathcal {C}}}

\newcommand{\CF}{{\mathcal {F}}}

\newcommand{\CH}{{\mathcal {H}}}

\newcommand{\CJ}{{\mathcal {J}}}

\newcommand{\CM}{{\mathcal {M}}}
\newcommand{\CN}{{\mathcal {N}}}

\newcommand{\CS}{{\mathcal {S}}}
\newcommand{\CT}{{\mathcal {T}}}

\newcommand{\CZ}{{\mathcal {Z}}}

\newcommand{\FH}{{\mathfrak {H}}}
\newcommand{\FI}{{\mathfrak {I}}}

\newcommand{\Fa}{{\mathfrak {a}}}

\newcommand{\Fc}{{\mathfrak {c}}}

\newcommand{\Fh}{{\mathfrak {h}}}

\newcommand{\RJ}{{\mathrm {J}}}

\newcommand{\RM}{{\mathrm {M}}}

\newcommand{\RO}{{\mathrm {O}}}

\newcommand{\RS}{{\mathrm {S}}}

\newcommand{\RW}{{\mathrm {W}}}

\newcommand{\Aut}{{\mathrm{Aut}}}

\newcommand{\End}{{\mathrm{End}}}
\newcommand{\el}{{\mathrm{ell}}}

\newcommand{\GL}{{\mathrm{GL}}}

\newcommand{\Ind}{{\mathrm{Ind}}}

\newcommand{\Id}{{\mathrm{Id}}}

\newcommand{\Lie}{{\mathrm{Lie}}}

\renewcommand{\Re}{{\mathrm{Re}}}

\newcommand{\Res}{{\mathrm{Res}}}

\newcommand{\SL}{{\mathrm{SL}}}

\newcommand{\Spec}{{\mathrm{Spec}}}

\newcommand{\Sym}{{\mathrm{Sym}}}

\newcommand{\Sp}{{\mathrm{Sp}}}

\newcommand{\st}{{\mathrm{st}}}

\newcommand{\tr}{{\mathrm{tr}}}

\newcommand{\ud}{\,\mathrm{d}}

\newcommand{\wh}{\widehat}

\newcommand{\bs}{\backslash}

\def\alp{{\alpha}}

\def\del{{\delta}}
\def\Del{{\Delta}}

\def\eps{{\epsilon}}

\def\sig{{\sigma}}

\def\std{\rm std}
\def\ome{{\omega}}
\def\Ome{{\Omega}}
\def\lam{{\lambda}}

\def\gam{{\gamma}}
\def\Gam{{\Gamma}}

\def\vpi{\varpi}

\newcommand{\Laf}{\operatorname{Laf}}
\newcommand{\LafSec}{\operatorname{LafSec}}
\newcommand{\fin}{\mathrm{fin}}
\newcommand{\triv}{\mathrm{triv}}

\def\vphi{\varphi}

\def\p{\prime}

\def\Tr{\mathrm{Tr}}

\newcommand{\sslash}{\mathbin{/\mkern-6mu/}}

\def\Irr{\mathrm{Irr}}
\def\der{\mathrm{der}}

\newtheorem{thm}{Theorem}[subsection]
\newtheorem{defin}[thm]{Definition}
\newtheorem{rmk}[thm]{Remark}

\newtheorem{pro}[thm]{Proposition}

\newtheorem{cor}[thm]{Corollary}

\newenvironment{num}{\begin{enumerate}}{\end{enumerate}}

\makeatletter

\newcommand{\Rmnum}[1]{\expandafter\@slowromancap\romannumeral #1@}
\makeatother

\newcommand*{\doublerightarrow}[2]{\mathrel{
  \settowidth{\@tempdima}{$\scriptstyle#1$}
  \settowidth{\@tempdimb}{$\scriptstyle#2$}
  \ifdim\@tempdimb>\@tempdima \@tempdima=\@tempdimb\fi
  \mathop{\vcenter{
    \offinterlineskip\ialign{\hbox to\dimexpr\@tempdima+1em{##}\cr
    \rightarrowfill\cr\noalign{\kern.5ex}
    \rightarrowfill\cr}}}\limits^{\!#1}_{\!#2}}}
\newcommand*{\triplerightarrow}[1]{\mathrel{
  \settowidth{\@tempdima}{$\scriptstyle#1$}
  \mathop{\vcenter{
    \offinterlineskip\ialign{\hbox to\dimexpr\@tempdima+1em{##}\cr
    \rightarrowfill\cr\noalign{\kern.5ex}
    \rightarrowfill\cr\noalign{\kern.5ex}
    \rightarrowfill\cr}}}\limits^{\!#1}}}
\makeatother

\def\bc{\mathbf{c}}

\def\temp{\mathrm{temp}}
\def\sep{\mathrm{sep}}
\def\heart{\heartsuit}

\firstpage{1}
\lastpage{23}

\title[
The Lafforgue Transform]{The Lafforgue Transform}

\author{Zhilin Luo}
\address{Department of Mathematics, Purdue University, West Lafayette, IN 47907 USA}
\email{luo642@purdue.edu}

\begin{document}

\begin{abstract}
In this report, we survey the \emph{Lafforgue transform}, a natural representation-theoretic map from the Kazhdan cocenter to the Bernstein center of a reductive group over a local field. Its motivation comes from our study of the Braverman-Kazhdan proposal, a conjectural framework for establishing the meromorphic continuation and functional equations of general automorphic $L$-functions through harmonic analysis on possibly singular varieties. For $G=\SL_2$ and $\GL_2$, we explain how the stable part of the Lafforgue transform and its distinguished section are realized as integral transforms on the Steinberg-Hitchin base, with formulas of motivic nature that are valid over every local field.
\end{abstract}

\maketitle

\tableofcontents

\section*{Introduction}\label{sec:intro}

One of the central themes in the Langlands program is the study of automorphic 
$L$-functions. These functions, which vastly generalize the Riemann zeta function,
encode arithmetic, geometric, and representation-theoretic information in a 
remarkably unified way. They are initially defined via Euler products and converge
only in a half-plane. One of the major implications of the Langlands functoriality
principle is that these $L$-functions should admit analytic continuation to the 
whole complex plane and satisfy functional equations \cite{langlandsproblems}.

A fundamental method for proving such analytic properties is harmonic analysis on
algebraic varieties. The classical prototype is Riemann's proof of the analytic continuation
and functional equation of the zeta function, where the Poisson summation formula for the
Gaussian plays the decisive role. This idea was systematically developed by J. Tate in his
thesis \cite{tatethesis}, and later by R. Godement and H. Jacquet for the standard
$L$-functions of general linear groups \cite{gjzeta} through the theory of local and 
global zeta integrals. Let $\RM_n$ be the variety of $n\times n$ matrices. It contains the group of invertible matrices $\GL_n$ as an open dense subvariety. For a global field $k$ with
adele ring $\BA_k$, set
\begin{equation}\label{eq:intro:GJ-zeta}
\CZ(s,\phi,\vphi_\pi)=
\int_{\GL_n(\BA_k)}
\phi(g)\vphi_\pi(g)|\det g|_{\BA}^{s+(n-1)/2}\ud g,
\end{equation}
where $\phi\in \CS(\RM_n(\BA_k))$ is a Schwartz-Bruhat function and $\vphi_\pi$ is a
matrix coefficient of a cuspidal automorphic representation $\pi$ of $\GL_n(\BA_k)$.
The classical Fourier transform $\CF_\psi$ on $\CS(\RM_n(\BA_k))$, where 
$\psi = \otimes_{\nu\in |k|}\psi_\nu:k\bs \BA_k\to \BC^\times$ is a nontrivial additive character, together with the Poisson
summation formula, implies the meromorphic continuation
and functional equation of \eqref{eq:intro:GJ-zeta}. The standard $L$-function
$L(s,\pi)$ is then recovered from these zeta integrals as poles of 
the family of tempered distributions $\{\CZ(s,\cdot,\vphi_\pi)\mid \vphi_\pi\in \CC(\pi)\}$
on $\RM_n(\BA_k)$.

Over the last several decades, various Rankin-Selberg integrals have been discovered, giving constructions of automorphic $L$-functions for special classes of
pairs $(G,\rho)$. Around 2000, A. Braverman and D. Kazhdan proposed a conjectural framework for such a theory for general $(G,\rho)$. Their proposal was
elaborated by L. Lafforgue \cite{lafforgue2016principe} and given in a more precise form by B.-C. Ng\^o \cite{ngo2016hankel}.
Roughly speaking, for a pair $(G,\rho)$, one expects a Schwartz space $\CS^\rho\big(G(\BA_k)\big)\simeq \bigotimes^\p_{\nu\in |k|}\CS^\rho\big(G(k_\nu)\big)$ together
with a $G(\BA_k)\times G(\BA_k)$-equivariant Fourier operator $\CF^\rho_\psi\simeq \bigotimes_{\nu\in |k|}\CF^\rho_{\psi_\nu}$. The expected Poisson summation
formula for the pair $\big(\CS^\rho\big(G(\BA_k)\big),\CF^\rho_\psi\big)$ should then imply the meromorphic continuation and functional equation of the general
automorphic $L$-functions attached to $(G,\rho)$.

One of the essential ingredients in the proposal is the Fourier operator $\CF^\rho_\psi$. Over a local field $k_\nu$ with $\nu\in |k|$, the Fourier operator $\CF^\rho_{\psi_\nu}$ is expected to be given by
convolution with a stably invariant distribution $\CJ^\rho_{\psi_\nu}$ on $G(k_\nu)$ that is represented by a locally integrable function smooth on the generic locus.
Up to mild modification, $\CJ^\rho_{\psi_\nu}$ is expected to lie in the stable Bernstein center of $G(k_\nu)$, in the sense that for any irreducible admissible representation $\pi_\nu$ of $G(k_\nu)$,
after regularization one has
$$
\pi_\nu(\CJ^\rho_{\psi_\nu}) = \gam(\cdot,\pi_\nu,\rho,\psi_\nu)\Id_{\pi_\nu}
$$
where the right-hand side is the Langlands local gamma factor. Because $\CJ^\rho_{\psi_\nu}$
is stably invariant, it descends to a function on the Steinberg-Hitchin base of $G$; in the case of $G=\SL_2$, this base is the
affine line with trace coordinate.

Our understanding of $\CJ^\rho_{\psi_\nu}$ is limited. For $(\GL_n,\rho=\std)$,
namely in Godement-Jacquet theory, $\CJ^\rho_{\psi_\nu} = \psi_\nu\big(\tr(\cdot)\big)$ is the standard
Fourier kernel. When $G$ is a torus, the kernel $\CJ^\rho_{\psi_\nu}$ is represented by
a smooth function on $G(k_\nu)$; this was mentioned in \cite{BK00}, made more precise in
\cite{ngo2016hankel}, and treated in detail in \cite{LN}. The author calculated the
spherical part of the distribution $\CJ^\rho_{\psi_\nu}$ for any $(G,\rho)$ over a local field 
of characteristic zero and showed that it can be inserted into the Arthur-Selberg trace formula \cite{MR3990815}.
In \cite{JLZ}, D. Jiang, L. Zhang, and the author give an explicit expression for the kernel for
$(\BG_m\times \Sp_{2n},\rho=\std)$. 
Nevertheless, for general $(G,\rho)$, the construction of $\CJ^\rho_{\psi_\nu}$ remains mysterious.

Motivated by Godement-Jacquet, Ng\^o speculated that the distribution $\CJ^\rho_{\psi_\nu}$
should be patched from the naive kernels $\{\CJ^\rho_{\psi_\nu,T}\}_{T\subset G}$ where 
$T$ runs over maximal tori of $G$ up to stable conjugation. However, as pointed out by 
Lafforgue \cite{lafforguegl2}, the naive patching construction fails to provide the correct kernel for $\rho=\Sym^2$ of $G=\BG_m\times \SL_2$
because it does not act on principal series representations with the correct gamma factors. For $G=\SL_2$ or $\GL_2$, Lafforgue proposed a recipe
to correct the naive patching construction so that the resulting distribution acts on the principal series representations by multiplication 
by the correct Langlands local gamma factor. However, whether the resulting distribution acts on supercuspidal representations by the correct gamma factors 
is unclear. Comparing with the result of \cite{JLZ}, one can show that Lafforgue's recipe works for $\rho=\Sym^2$ of $G=\BG_m\times \SL_2$, but it fails for supercuspidal representations and
higher symmetric powers. 

In \cite{LN}, we develop a general framework for the construction of
$\CJ^\rho_{\psi_\nu}$, based on Langlands stable transfer \cite{MR3117742} and the 
so-called \emph{Lafforgue transform}. These transforms exist for general reductive algebraic groups and can be realized as explicit analytic formulas
on the Steinberg-Hitchin base for $G=\SL_2$ and $\GL_2$.

Among these ingredients, the Lafforgue transform is the central object of this report. It is a 
natural representation-theoretic map 
\begin{equation}\label{eq:intro:Laf}
\Laf:\CC_{G(k_\nu)}\longrightarrow \CZ_{G(k_\nu)}
\end{equation}
from the Kazhdan cocenter to the Bernstein center of $G(k_\nu)$. Here the Bernstein center 
can be viewed as the algebra of multipliers on smooth representations of $G(k_\nu)$, while 
the Kazhdan cocenter can be viewed as the completion of the usual cocenter of the Hecke algebra, or
equivalently as the space of orbital integrals of the Hecke algebra. Concretely, for $c\in \CC_{G(k_\nu)}$ and
$[M(k_\nu),\sig_\nu]$ a cuspidal datum, the Lafforgue transform is characterized by the 
formula 
\begin{equation}\label{eq:intro:Laf-spectral}
\Laf(c)\big([M(k_\nu),\sigma_\nu]\big)=
\Tr\big(\Ind(\sig_\nu),c\big),
\end{equation}
where the right-hand side is the trace evaluation of the orbital integral represented by $c\in \CC_{G(k_\nu)}$
on the fully induced representation $\Ind(\sig_\nu)$.

For $G=\SL_2$ or $\GL_2$, this representation-theoretic map becomes a concrete 
analytic formula. After restricting to the stable finite part, there is a 
distinguished section 
$$
\LafSec:\CZ_{G(k_\nu)}^{\st,\fin}\to \CC_{G(k_\nu)}^{\st,\fin}
$$
of the Lafforgue transform that can be realized as an explicit integral transform 
on the Steinberg-Hitchin base. For $z\in \CZ^{\st,\fin}_{G(k_\nu)}$ that is represented
by a stable function $\CJ_{G(k_\nu)} = \CJ_{\Fc_\nu}\circ \Fc_\nu$, where $\Fc_\nu:G(k_\nu)\to \Fc(k_\nu)$
is the Chevalley quotient, given by the trace for $G=\SL_2$ and by the trace and determinant for $G=\GL_2$,
the corresponding element in the Kazhdan cocenter is represented by
$$
\LafSec(z) = 
\CF_{\psi_\nu}\bigg(
\frac{\CF_{\psi^{-1}_\nu}(\CJ_{\Fc_\nu})}{|\cdot|_\nu}  
\bigg)
$$
where the Fourier transforms are taken in the trace variable of the Steinberg-Hitchin base. 
This formula is one of the main points of \cite{LN}. It shows that, although the
Lafforgue transform is defined representation-theoretically, its stable section in the rank-one
case is governed by an explicit analytic formula that is of motivic nature and valid
over every local field, including fields of residual characteristic two. For general reductive groups,
Ng\^o and the author have formulated conjectural analytic formulas for the stable
Lafforgue transform and its section on the Steinberg-Hitchin base; see the work in
preparation \cite{NL-general}.

The other ingredient is Langlands stable transfer, which provides a canonical map from the Kazhdan cocenter of $G(k_\nu)$ to the Kazhdan cocenters of its maximal tori.
Again in the rank-one case,
Langlands stable transfer can be realized by a simple analytic formula in the trace variable, closely connected with the Gelfand-Graev character formula \cite{MR0155931}. In particular, after
applying the composition of $\LafSec$ with Langlands stable transfer to the conjectural Fourier kernel $\CJ^{\rho}_{\psi_\nu}$,
one obtains the naive kernel predicted by Ng\^o's patching construction.

The purpose of this report is to introduce the Lafforgue transform. We first recall the Bernstein center, the Kazhdan cocenter, and their stable variants. We then define the Lafforgue transform and explain its basic properties. Finally, in the cases $G=\SL_2$ and $G=\GL_2$, we explain how its stable part and distinguished section can be written as explicit integral transforms on the Steinberg-Hitchin base.

\subsection*{Acknowledgements}

The author would like to thank the organizers of the International
Congress of Chinese Mathematicians for inviting him to speak at the conference.

\section{Bernstein center}\label{sec:BC}

In this section, we review the definition of the Bernstein center, together with its
archimedean and stable variants.

\subsection{Non-archimedean case}\label{subsec:non-archi:BC}

Let $G$ be a reductive group over a non-archimedean local field $F$. Following \cite{BD84}, the Bernstein center $\CZ=\CZ_G$ admits the following equivalent definitions:

\begin{num}
\item\label{num:BC:1} It is the ring of endomorphisms of the identity functor on the category of smooth representations of $G(F)$;

\item\label{num:BC:2} It is the projective limit of the center of the Hecke algebra at finite levels $\CH_K = e_K\CH e_K$ where $e_K$ runs over the characteristic measures of open compact subgroups $K$ of $G(F)$;

\item\label{num:BC:3} It is the ring of invariant distributions $\CJ$ on $G(F)$ that are essentially compactly supported, i.e. $\CJ e_K = e_K\CJ\in \CH$ for any open compact subgroup $K$ of $G(F)$;

\item\label{num:BC:4} It is the ring of regular functions on the Bernstein variety $\Ome = \Ome_G$, whose points parametrize $G(F)$-conjugacy classes of pairs $(M,\sig)$, where $M$ is a Levi subgroup of $G$ and $\sig$ is an irreducible supercuspidal representation of $M(F)$.
\end{num}

The Bernstein variety $\Ome$ is a disjoint union of infinitely many finite-dimensional algebraic varieties. It has infinitely many connected components parametrized by the inertial equivalence classes of such pairs. Let $\CZ^\fin$ be the nonunital subalgebra of $\CZ$ whose elements are regular functions on $\Ome$ that are supported on finitely many connected components. It is also an ideal in $\CZ$.

For every finite subset $I$ of $\pi_0(\Ome)$, where $\pi_0(\Ome)$ denotes the set of connected components of $\Ome$, let $\eps_I\in \CZ$ be the regular function on $\Ome$ supported on $\bigcup_{i\in I}\Ome_i$ and taking the value one on those connected components. We abuse the notation and still write $\eps_I$ for the corresponding distribution on $G(F)$. For two finite subsets $I\subset J$ of $\pi_0(\Ome)$, the identity $\eps_I = \eps_I\eps_J$ induces a homomorphism of unital algebras $\CZ\eps_J\to \CZ\eps_I$ given by $\CJ\mapsto \CJ\eps_I$. There is also a homomorphism of nonunital algebras $\CZ\eps_I\to \CZ\eps_J$. The Bernstein center can also be represented as the projective limit of $\CZ\eps_I$, and $\CZ^\fin$ can be represented as the injective limit of $\CZ\eps_I$. Given a distribution $\CJ$ in the Bernstein center $\CZ$, let $\gam(\CJ,[\pi])$ denote the value of the corresponding regular function on $\Ome$ at the supercuspidal support $[\pi]$ of $\pi$. For the following function-valued formula, assume that $\CJ\in\CZ^\fin$, so that its spectral multiplier is supported on finitely many Bernstein components. The relation between $\CJ$ and the function $\gam(\CJ,[\cdot])$ is then expressed by the Plancherel formula
\begin{align}\label{eq:BC:1}
\CJ(g) = \int_{\pi\in \Pi_G^\temp}
\gam(\CJ,[\pi]) \Theta_\pi(g)\ud \pi
\end{align}
where $\Pi^\temp_G$ is the tempered dual of $G(F)$, $\ud\pi$ is the Plancherel measure, and $\Theta_\pi$ is the distribution character of $\pi$. For a general element $\CJ\in\CZ$, Plancherel inversion is understood distributionally, after pairing against test functions.

\subsection{Archimedean variant}\label{subsubsec:archiBC}

Let $G$ be a reductive algebraic group over an archimedean local field $F$. By restriction of scalars, we may take $F=\BR$. We work with the category of Casselman-Wallach representations, i.e. admissible smooth Fréchet representations of moderate growth. By \cite[Prop.~2.20]{bkglobalization}, this category is equivalent to the category of continuous non-degenerate $\CS$-modules, where $\CS = \CS\big(G(F)\big)$ is the Schwartz algebra of rapidly decreasing functions on $G(F)$.

A natural candidate for the archimedean analogue of the Bernstein center is the algebra of multipliers of $\CS$, i.e. the ring of continuous endomorphisms $\End_{\CS\times \CS}(\CS)$, which is also the algebra of invariant linear functionals on $G(F)$ sending $\CS$ to $\CS$ via convolution. A spectral characterization of $\End_{\CS\times \CS}(\CS)$ relies on a trace Paley-Wiener theorem for $\CS$; in the form needed here, this seems to be known only for $G=\SL_2$ (see \cite[Remark~3.4.1]{LN} and the discussion above for a summary). Following \cite[\S 2]{MR4298750}, we introduce a subalgebra of $\End_{\CS\times \CS}(\CS)$ using a spectral characterization parallel to \eqref{eq:BC:1}, which is sufficient for many purposes.

Let $G_\BC = G\otimes_\BR\BC$. We introduce the following datum associated to $G_\BC$.

\begin{itemize}
\item Let $(\BX^*,\Phi,\BX_*,\Phi^\vee)$ be the root datum for $G_\BC$, i.e. $\BX^*$ and $\BX_*$ are the weight and coweight lattices of the abstract Cartan subgroup of $G_\BC$, with $\Phi$ and $\Phi^\vee$ the sets of roots and coroots of $G_\BC$, respectively;

\item Let $\RW$ be the Weyl group of $(\BX^*,\Phi,\BX_*,\Phi^\vee)$;

\item Let
\[
\mathscr A_G=\Aut(\BX^*,\Phi,\BX_*,\Phi^\vee)^\heart
\]
be the finite set of automorphisms of $(\BX^*,\Phi,\BX_*,\Phi^\vee)$ of order at most $2$; it is stable under the action of $\RW$;

\item Let $\Fh^* = \BX^*\otimes_\BZ \BR$;

\item For $\vartheta\in \mathscr A_G$, let $\Fh^*_\vartheta$ (resp. $\Fh^{*-}_\vartheta$) be the $+1$-eigenspace (resp. $-1$-eigenspace) of the action of $\vartheta$ on $\Fh^*$. Let $\BY^*_\vartheta$ be the projection of $\rho+\BX^*\subset \Fh^*$ onto $\Fh^{*-}_\vartheta$, where $\rho$ is the half sum of positive roots in $\Phi$ with respect to a choice of base of $\Phi$.
\end{itemize}
$\BY^*_\vartheta$ is a translation of a discrete subgroup of $\Fh^{*-}_\vartheta$. Although $\BY^*_\vartheta$ is not necessarily a subgroup, below we sometimes write $\BH\oplus \BY^*_\vartheta$ for the subset of $\Fh^*_\BC$ consisting of elements of the form $\alp+\vpi$, where $\alp\in \BH$ and $\vpi\in \BY^*_\vartheta$.

Following \cite[\S 2]{MR4298750}, for a subset $\theta\subset \mathscr A_G$ that is stable under the conjugation action of $\RW$, we introduce the following spaces of holomorphic functions on $\Fh^*_\BC$.

\begin{defin}\label{defin:Schwartzmultiplier}
Consider the following spaces of holomorphic functions on $\Fh^*_\BC$.
\begin{enumerate}
\item Let $\CM_\theta(\Fh^*_\BC)$ be the space of holomorphic functions $\mu$ on $\Fh^*_\BC$ such that for every $\vartheta\in \theta$ and every $M>0$, there exists $r_{\vartheta,M}\in \BR$ such that 
\begin{equation}\label{eq:schwartzmultiplier:1}
\sup_{
\substack{
\alp\in \Fh^*_{\vartheta,\BC},||\Re\alp||<M
\\
\vpi\in \BY^*_\vartheta
}
}
|\mu(\alp+\vpi)|
\cdot 
(1+||\alp+\vpi||)^{r_{\vartheta,M}}<\infty
\end{equation}
holds for some, hence every, Euclidean norm $||\cdot||$ on $\Fh^*_\BC$;

\item Let $\CN_\theta(\Fh^*_\BC)$ be the space of holomorphic functions $\mu$ on $\Fh^*_\BC$ such that for every $\vartheta\in \theta$, every $M>0$, and every $r\in \BR$, 
\begin{equation}\label{eq:schwartzmultiplier:2}
\sup_{
\substack{
\alp\in \Fh^*_{\vartheta,\BC},
||\Re\alp||<M
\\
\vpi\in \BY^*_\vartheta
}
}
|\mu(\alp+\vpi)|
\cdot 
(1+||\alp+\vpi||)^r<\infty
\end{equation}
holds for some, hence every, Euclidean norm $||\cdot||$ on $\Fh^*_\BC$;

\item Define $\CM^\natural_\theta(\Fh^*_\BC)$ to be the subspace of $\CM_{\theta\cup\{1\}}(\Fh^*_\BC)$ consisting of functions $\mu$ such that, for every $\vartheta\in \theta$, $\mu|_{\Fh^*_{\vartheta,\BC}+\vpi} = 0$ for all but finitely many elements $\vpi\in \BY^*_\vartheta$.
\end{enumerate}
\end{defin}
When $\theta = \{1\}$, we suppress it in the subscripts in the above notation.

\begin{rmk}
We make the following comments on the definition.

\begin{enumerate}
\item In the definition of $\CM_\theta(\Fh^*_\BC)$ (resp. $\CN_\theta(\Fh^*_\BC)$), for every $\vartheta\in \theta$ and every $\vpi\in \BY^*_\vartheta$, the function $\alp\mapsto \mu(\alp+\vpi)$ has moderate vertical growth (resp. rapid decay on vertical strips) on $\Fh^*_{\vartheta,\BC}$. If the growth condition is restricted to $i\Fh^*_\vartheta\oplus \BY^*_\vartheta$, then it says that $\mu$ has polynomial growth;

\item The spaces $\CM_\theta(\Fh^*_\BC)$, $\CN_\theta(\Fh^*_\BC)$, and $\CM^\natural_\theta(\Fh^*_\BC)$ are all closed under multiplication and the action of $\RW$;

\item If a holomorphic function $\mu$ on $\Fh^*_\BC$ has moderate vertical growth, then $\mu\cdot \CM^\natural_\theta(\Fh^*_\BC)\subset \CM^\natural_\theta(\Fh^*_\BC)$. In particular, we have $\CM(\Fh^*_\BC)\cdot \CM^\natural_\theta(\Fh^*_\BC)\subset \CM^\natural_\theta(\Fh^*_\BC)$. Note that $\CM(\Fh^*_\BC) = \CM_{\{1\}}(\Fh^*_\BC)$.
\end{enumerate}
\end{rmk}

Let $\theta_G$ be the subset of $\mathscr A_G$ consisting of elements of the form $w\vartheta$, where $w\in \RW$ and $\vartheta\in \mathscr A_G$ induces the real form $G$ of $G_\BC$. Then $\theta_G$ is stable under the action of $\RW$.

We introduce the following definitions.

\begin{defin}\label{defin:archimedean-Bernstein-center}
Let $\CZ = \CM_{\theta_G\cup\{1\}}(\Fh^*_\BC)^\RW$, and $\CZ^\fin = \CM^\natural_{\theta_G\cup\{1\}}(\Fh^*_\BC)^\RW$.
\end{defin}
$\CZ$ will play the role of the Bernstein center in the archimedean setting.

Following \cite[Thm.~2.13]{MR4298750}, any element of $\CZ$ gives rise to a multiplier in $\End_{\CS\times \CS}(\CS)$. Moreover, this assignment gives a homomorphism of complex algebras $\CZ\to \End_{\CS\times \CS}(\CS)$. By the matricial Paley-Wiener theorem for smooth and compactly supported functions on $G(F)$ \cite{MR2183287}, any multiplier $\gam\in \CZ$ uniquely determines an invariant distribution $\CJ_\gam$ on $G(F)$, giving the analogue of \eqref{eq:BC:1} in the archimedean setting.

\begin{rmk}\label{rmk:defin:Schwartz-multiplier}
In \cite{LN}, we define $\CZ = \CN_{\theta_G\cup\{1\}}(\Fh^*_\BC)^\RW$, which is sufficient for the purposes of \emph{loc. cit.} In this report, we work with $\CM_{\theta_G\cup \{1\}}(\Fh^*_\BC)^\RW$ for general considerations.
\end{rmk}

\subsection{Stable Bernstein center}\label{subsec:SBC}

Let $G$ be a reductive group over a local field $F$. Following \cite{bezrukavnikov2013categorical}, let $\CZ^\st$ be the subset of $\CZ$ consisting of stable distributions on $G(F)$ in the sense of \cite{MR0540901}, where an invariant distribution is called stable if for any test function $h\in \CH$, $\CJ(h)= 0$ whenever the stable orbital integral of $h$ vanishes identically. Let $\CZ^{\st,\fin} = \CZ^\st\cap \CZ^\fin$. The stable center conjecture of \cite{bezrukavnikov2013categorical} stipulates that $\CZ^\st$ is a subalgebra of $\CZ$, at least when $F$ is non-archimedean.

There are several related definitions of the stable Bernstein center. We recall them below.

\begin{num}
\item Let $\FI_\st\subset \CH$ be the subspace of those test functions whose stable orbital integrals vanish identically. Following \cite{MR2983012}, we define $\CZ^{\st,\p}$ to be the subalgebra of $\CZ$ consisting of distributions that stabilize $\FI_\st$, i.e.
$$
\CZ^{\st,\p} = \{z\in \CZ\mid z*\FI_\st\subset \FI_\st\}.
$$
Comparing with the definition of $\CZ^\st$, the definition of $\CZ^{\st,\p}$ has the benefit that it is immediately a subalgebra of $\CZ$.

\item In general, as shown in \cite[Prop.~3.5.3]{LN}, $\CZ^{\st,\p}\subset \CZ^\st$ and when $G=\SL_2$ over a non-archimedean local field of odd residual characteristic, $\CZ^\st = \CZ^{\st,\p}$.

\item Let $\Ome^\st = \Ome^\st_G$ be the variety of infinitesimal characters of $G$ \cite[\S 5]{MR3444233}\cite[\S 6]{MR2983012}. Assuming that the local Langlands conjecture is known for $G$ and its Levi subgroups, there exists a finite-to-one map $p:\Ome\to\Ome^\st$ sending the cuspidal datum in $\Ome$ to the corresponding $\wh{G}$-conjugacy classes of semisimple local $L$-parameters. We define $\CZ^{\st,\p\p}$ to be the subalgebra of $\CZ$ consisting of regular functions on $\Ome$ that are the pullback of regular functions on $\Ome^\st$ along $p$.
\end{num}

In general, the relation among $\CZ^\st$, $\CZ^{\st,\p}$, and $\CZ^{\st,\p\p}$ is not clear. When $G=\SL_2$ over a non-archimedean local field of odd residual characteristic, one can show that the three notions coincide \cite[Prop.~3.5.3]{LN}.

\section{Kazhdan cocenter}\label{sec:KCC}

In this section, we review the notion of the Kazhdan cocenter introduced in \cite[\S 3]{LN}.

\subsection{Normalization of measures}\label{subsec:measures}

We follow \cite{FLN10} and \cite[\S 3.1]{LN} to define measures on the spaces of (stable) conjugacy classes of $G(F)$.

We first consider the case in which $G$ is split, semisimple, and simply connected, as treated in \cite[\S 3.1]{LN}. By \cite{Steinberg-conjugacy}, the invariant-theoretic quotient of $G$ acting on itself via conjugation is given as follows
$$
\Fc=\Fc_G=G\sslash G=\Spec(F[G]^G) = \BA^r
$$
with coordinates $c_1,\ldots,c_r$ given by $c_i(g) = \tr(\rho_{\ome_i}(g))$, where the $\rho_{\ome_i}$ are the fundamental representations of $G$. Following \cite{FLN10}, we call $\Fc$ the Steinberg-Hitchin base of $G$. By the Chevalley restriction theorem, there is an isomorphism
$$
\Fc=T\sslash \RW = \Spec(F[T]^\RW)
$$
where $T$ is a maximal split torus and $\RW$ is the corresponding Weyl group. We define the Weyl discriminant $D\in F[\Fc]$ as the $\RW$-invariant function on $T$ given by 
$$
D(t) = \prod_{\alp\in \Phi}(1-\alp(t))
$$
where $\Phi = \Phi(G,T)$ is the set of roots. The absolute value $|D(t)|$ of the Weyl discriminant $D(t)$ admits a square-root constructed as follows. Choose a $\BZ$-basis $\chi_1,\ldots,\chi_r$ of the group $\BX^*(T)$ of rational characters of $T$. We then have an invariant $r$-form on $T$ given by
\begin{align}\label{eq:normalize:1}
\frac{\ud \chi_1}{\chi_1}
\wedge\cdots\wedge\frac{\ud \chi_r}{\chi_r}
\end{align}
which, up to a sign, is independent of the choice of the $\BZ$-basis $\chi_1,\ldots,\chi_r$ of $\BX^*(T)$. By pulling back the $r$-form $\ud c_1\wedge\cdots\wedge \ud c_r$ along the morphism $\pi_T:T\to \Fc = T\sslash \RW$, we have
\begin{align}\label{eq:normalize:2}
\pi^*_T
\ud c_1\wedge\cdots\wedge \ud c_r =
\Del_T(t)
\frac{\ud \chi_1}{\chi_1}
\wedge\cdots\wedge
\frac{\ud \chi_r}{\chi_r}
\end{align}
where $\Del_T(t)$ is an algebraic function on $T$ well-defined up to a sign. Its absolute value 
$$
\Del(t) = |\Del_T(t)|
$$
is a real-valued function invariant under the action of the Weyl group $\RW$. By \cite[Prop.~3.29]{FLN10},
$$
|D(t)| = \Del(t)^2.
$$
The $r$-form $\ud c_1\wedge\cdots\wedge \ud c_r$ on the $r$-dimensional affine space $\Fc$ determines a measure
\begin{align}\label{eq:normalize:3}
\ome_\Fc = |\ud c_1\wedge\cdots\wedge \ud c_r|
\end{align}
on $\Fc(F) = F^r$. It follows from \eqref{eq:normalize:2} that $\Del^{-1}\pi^*_T\ome_\Fc$ is an invariant measure on $T(F)$ that coincides with the measure associated with the invariant $r$-form \eqref{eq:normalize:1}. If $T$ is a maximal torus of $G$ but not necessarily split, the finite map $\pi_T:T\to \Fc$ remains available, and the formula $\Del^{-1}\pi^*_T\ome_\Fc$ again defines an invariant measure on $T(F)$. We denote the corresponding measure by
$$
\ome_T = \Del^{-1}\pi^*_T\ome_\Fc.
$$

Next we consider the quasi-split, semisimple, simply connected case.
A quasi-split form of a reductive group $G$ is determined by a $1$-cocycle $\vphi:\sig\in \Gam_F\mapsto \vphi(\sig)$ valued in the automorphism group of the Dynkin diagram of $G$, where $\Gam_F$ is the absolute Galois group of $F$. As a result, $\Gam_F$ acts on the set of fundamental representations
$$
\sig:\rho_{\ome_i}
\mapsto \rho_{\sig(\ome_i)},\quad \rho_{\sig(\ome_i)}(g) = 
\rho_{\ome_i}(\vphi(\sig)g),
$$
and consequently on the trace characters $c_i = \tr(\rho_{\ome_i})$. Observe that
$$
\rho_{\sig(\tau(\ome_i))}(g) = 
\rho_{\tau(\ome_i)}(\vphi(\sig)g) = 
\rho_{\ome_i}(\vphi(\tau)\vphi(\sig)g) = 
\rho_{\ome_i}(\vphi(\tau\sig)g) = 
\rho_{(\tau\sig)(\ome_i)}(g).
$$
Let $G_\vphi$ be the quasi-split group determined by the cocycle $\vphi:\sig\in \Gam_F\mapsto \vphi(\sig)$. The action of $\Gam_F$ on $G_\vphi(\bar{F}^\sep)$, which as a set is just $G(\bar{F}^\sep)$, is given by 
$$
\sig_\vphi:g\mapsto \vphi(\sig)(\sig(g)) = \sig(\vphi(\sig)g).
$$
Here $\bar{F}^\sep$ is the separable closure of $F$. Suppose that $\{h_\sig\}$ is any cocycle, with respect to this action, valued in the center of $G_\vphi(\bar{F}^\sep) = G(\bar{F}^\sep)$. Then $\rho_{\ome_i}(h_\sig)\in \bar{F}^\sep$ acts as a scalar $\eps_i(\sig)$. Moreover, $\sig(\eps_{\sig(\ome_i)}(\tau))\eps_i(\sig) = \eps_i(\tau\sig)$, because 
\begin{align*}
\sig(\rho_{\sig(\ome_i)}(h_\tau))\rho_{\ome_i}(h_\sig)
&=
\rho_{\sig(\ome_i)}(\sig(h_\tau))\rho_{\ome_i}(h_\sig) \\
&=
\rho_{\ome_i}(\vphi(\sig)(\sig(h_\tau))h_\sig) \\
&= \rho_{\ome_i}(h_{\tau\sig}).
\end{align*}
As a consequence, we twist the linear space over $F$ for the coordinates $c_1,\ldots,c_r$ using the cocycle
$$
\sig:(c_1,\ldots,c_r)\mapsto (c_1^\p,\ldots,c_r^\p),\quad c^\p_i = \eps_i(\sig)c_{\sig(\ome_i)}.
$$
Taking $\{h_\sig\}$ to be the trivial cocycle, for any inner form of the quasi-split form $G_\vphi$, we obtain a twisted form of the Steinberg-Hitchin base of the semisimple, simply connected, split outer form of $G$ that is a vector space over $F$ with a measure inherited from \eqref{eq:normalize:3}.

In general, by \cite{kott2}, any reductive group $G$ admits a $z$-extension $\alp:G_z\to G$ with the following properties:
\begin{itemize}
\item $G_z$ is a connected reductive group over $F$ with simply connected derived group;

\item $\alp$ is surjective;

\item $\ker \alp$ is central in $G_z$ and is isomorphic to a product of the form $\Res_{L/F}\BG_m$ (where each $L$ is a finite extension of $F$ splitting $G$).
\end{itemize}
By pulling back along $\alp$, we reduce the problem to the case in which $G_\der$ is simply connected. Let $Z$ be the identity component of the center of $G$. There is an exact sequence
\begin{equation}\label{eq:normalize:4}
1\to A\to Z\times G_\der\to G\to 1,\quad A=G_\der \cap Z.
\end{equation}
Following \cite{FLN10}, to simplify the exposition we assume that $A$ is étale over $F$. By \cite[Prop.~2.2]{MR207712}, for $g\in G(\bar{F}^\sep)$ whose image lies in $G(F)$, we can write it as a product $g=g_1g_2$ with $g_1\in G_\der(\bar{F}^\sep)$ and $g_2\in Z(\bar{F}^\sep)$. Since $g\in G(F)$, for any $\sig\in \Gam_F$, we have $\sig(g_1)^{-1}g_1 = \sig(g_2)g_2^{-1}\in A(\bar{F}^\sep)$. Since $g_1$ is uniquely determined up to an element of $A(\bar{F}^\sep)$, only its image $\eta(g)$ of the cocycle $\{\sig(g_1)^{-1}g_1\}$ in $H^1(F,A)$ is defined. Let $\FH$ be the subset of $H^1(F,A)$ obtained in this way. The Steinberg-Hitchin base of $G$ is defined to be
\begin{align}\label{eq:normalize:5}
\Fc = \bigsqcup_{\eta\in \FH}
\Fc_{G_\der,\eta}(F)\times Z_\eta(F)/A(F)
\end{align}
where $\Fc_{G_\der,\eta}$ is a twisted linear space of the underlying Steinberg-Hitchin base of $G_\der$ on which $A(F)$ acts as a subgroup of the center of $G_\der(F)$, and $Z_\eta$ is a $Z$-torsor defined by $\eta$ on which $A(F)$ acts as a subgroup of $Z(F)$. The space \eqref{eq:normalize:5} is also equipped with a natural measure, where the measure on $\Fc_{G_\der,\eta}(F)$ is inherited from \eqref{eq:normalize:3}, the measure on $Z_\eta(F)$ is the measure on the split torus, and the measure on $\Fc_{G_\der,\eta}(F)\times Z_\eta(F)/A(F)$ is given the quotient measure.

For $G$ semisimple and simply connected, let $\Fc^\p$ be the complement of the divisor defined by the Weyl discriminant. Its inverse image $G^\p$ in $G$ consists of (strongly) regular semisimple elements of $G$. For a general reductive group $G$, let $\Fc^\p$ be 
$$
\bigsqcup_{\eta\in \FH}
\Fc^\p_{G_\der,\eta}(F)\times Z_\eta(F)/A(F)
$$
The inverse image $G^\p$ in $G$ consists of strongly regular semisimple elements of $G$. The map 
\begin{equation}\label{eq:normalize:6}
\bc:G^\p(F)\to \Fc^\p(F)
\end{equation}
has nonempty fibers given by stable conjugacy classes. If $G$ is quasi-split, it is surjective by \cite[Thm.~4.1]{kott2}, using the preceding reduction to a simply connected derived group; splitness is not required. For a general $G$, surjectivity need not hold. In that case, in the following discussion of conjugacy classes we use $\Fc^\p(F)$ to denote the open image $\bc(G^\p(F))$. Thus $\Fc^\p(F)$ can be identified with the set of stable strongly regular semisimple conjugacy classes. This set carries a canonical structure of an $F$-analytic variety induced from $\Fc(F)$, together with the measure induced from $\ome_\Fc$. Let $C^\p(F)$ denote the set of strongly regular semisimple conjugacy classes in $G^\p(F)$. The fibers of the quotient map
\begin{equation} \label{eq:normalize:7}
    \mathbf{d}:G^\p(F) \to C^\p(F)
\end{equation}
are strongly regular semisimple conjugacy classes. We consider the map 
$$
\mathbf{e}:C^\p(F)\to \Fc^\p(F)
$$
from the set of strongly regular semisimple conjugacy classes of $G^\p(F)$ into the set of stable strongly regular semisimple conjugacy classes. This map is finite-to-one, since every stable strongly regular semisimple conjugacy class contains only finitely many conjugacy classes. It is also a covering: for every $x\in C^\p(F)$, there exist open neighborhoods $U$ of $x$ and $V$ of $\mathbf{e}(x)$ such that $\mathbf{e}$ induces a homeomorphism $U\to V$. The canonical measure $\ome_{\Fc}$ on $\Fc^\p(F)$ gives rise to a canonical measure $\ome_C$ on $C^\p(F)$. We also note that the map \eqref{eq:normalize:7} is an $F$-analytic submersion. Since $\Fc^\p(F)$ is equipped with a structure of an $F$-analytic manifold, its covering $C^\p(F)$ is also an $F$-analytic manifold. The measure $\ome_\Fc$ on $\Fc(F)$ pulls back to a measure on $C^\p(F)$. This extends to a measure on the set $C(F)$ of conjugacy classes of $G(F)$ by declaring that the complement of $C^{\p}(F)$ in $C(F)$ has measure zero.

\subsection{Normalization of stable orbital integrals and orbital integrals}\label{subsec:normalizationSOI-OI}

Let $\CH$ be the algebra of smooth and compactly supported functions on $G(F)$, where $G$ is a reductive group over $F$. When $F$ is archimedean, it is also convenient to work with $\CS$, the Schwartz algebra of rapidly decreasing functions on $G(F)$. In what follows, we focus on functions in $\CH$, and the discussion applies to functions in $\CS$ when $F$ is archimedean.

For every $h\in \CH$, consider the push-forward measure $\bc_!(h\ud g)$ defined so that for every continuous function $f:\Fc(F)\to \BC$ the following identity holds; here $\ud g$ is a Haar measure on $G(F)$:
$$
\langle \bc_!(h\ud g),f\rangle = \langle h\ud g,\bc^*f\rangle.
$$
Here $\langle \cdot,\cdot \rangle$ is the standard pairing between compactly supported measures and continuous functions and $\bc^*f = f\circ \bc$ is the pullback of $f$ along $\bc$. Based on the work of Harish-Chandra, there exists a function $\RS\RO(h)$ locally integrable with respect to $\ome_\Fc$ such that
$$
\bc_!(h\ud g) = \RS\RO(h)\ome_\Fc.
$$
When restricted to the strongly regular semisimple locus, 
$$
\RS\RO(h):\Fc^\p(F)\to \BC
$$
is a smooth function that is called the stable orbital integral of $h$. This is the normalization of stable orbital integrals adopted in \cite[(3.31)]{FLN10}.

We can apply the same construction for the map $G^\p(F)\to C^\p(F)$ and thus define a function 
$$
\RO(h):C^\p(F)\to \BC,
$$
such that 
$$
\mathbf{d}_!(h\ud g) = \RO(h)\ome_C,
$$
and we call it the orbital integral of $h$. This is equivalent to requiring that for every continuous function $f:C(F)\to \BC$ we have 
$$
\langle 
\mathbf{d}_!(h\ud g),f\rangle = 
\langle h\ud g,\mathbf{d}^*f\rangle.
$$
Now let us compare the above normalization of orbital and stable orbital integrals with the more usual definition. In the customary setting, we define the orbital integral for an element $\gam\in G^\p(F)$ with the formula 
\begin{equation}\label{eq:normalize-SOI:1}
\RO_\gam(h) = 
\int_{G(F)/G_\gam(F)}
h(g\gam g^{-1})\frac{\ud g}{\ud t}
\end{equation}
where $\ud t$ is a Haar measure of $T(F) = G_\gam(F)$. Let $T^\p\subset T$ be the strongly regular semisimple locus. By restricting to $T$, we have a function
\begin{equation}\label{eq:normalize-SOI:2}
\RO_T(h):T^\p(F)\to \BC
\end{equation}
such that for any $\gam\in T^\p(F)$, $\RO_T(h)(\gam) = \RO_\gam(h)$. Notice that there is a natural projection $p_T:T^\p(F)\to C^\p(F)$ that is a covering map. From \cite[\S 3.2]{LN}, the functions $\{\RO_T(h)\}_{T\in \CT}$ can be glued together to get a function $\RO^\CT(h)$ on $C^\p(F)$ where $\CT$ is the set of $G(F)$-conjugacy classes of maximal tori of $G(F)$. The following statement is proved in \cite[Cor.~3.2.2]{LN} using the Weyl integration formula.

\begin{cor}\label{cor:normalize:SOI}
We have the identity $\RO(h)(a) = \Del(a)\RO^\CT(h)(a)$ for $a\in C^\p(F)$.
\end{cor}

\subsection{Cocenter and orbital integrals}\label{subsec:CC-OI}

Following the work of Harish-Chandra \cite{HC70}\cite{varrealHC}, the function 
$$
\RO(h):C^\p(F)\to \BC
$$
is bounded for every $h\in \CH$. When $F$ is non-archimedean, the theorem of Kazhdan \cite{KazCus} shows that the kernel of the map $h\in \CH\mapsto \RO(h)$ is given by the commutator $[\CH,\CH]$. This induces an isomorphism of the cocenter $\CH/[\CH,\CH]$ with the image $\RO(\CH) = \{\RO(h)\mid h\in \CH\}$ when $F$ is non-archimedean. Similarly, the map $h\in \CH\mapsto \RS\RO(h)$ defines a map from the space of smooth and compactly supported functions to the space of stable orbital integrals, which we denote by $\RS\RO(\CH)$. We call $\RS\RO(\CH)$ the stable cocenter of $\CH$. When $F$ is archimedean, we also introduce the notation $\RO(\CS)$ (resp. $\RS\RO(\CS)$) for the space of orbital integrals (resp. stable orbital integrals) of functions in the Schwartz algebra $\CS$.

For a smooth irreducible admissible representation $\pi$ of $G(F)$, its distribution character $\Theta_\pi$ is, by definition, the distribution on $G(F)$ such that
$$
\tr(\pi(h)) = \int_{G(F)}
h(g)\Theta_\pi(g),\quad h\in \CH.
$$
The trace distribution extends to $\CS$ when $F$ is archimedean \cite[\S 8.1]{Wallach-RealI}.
From the work of Harish-Chandra, it is known that $\Theta_\pi$ is an invariant locally integrable function smooth on the regular semisimple locus. We abuse the notation and write $\Theta_\pi$ for the corresponding invariant function on $G(F)$.

Let $\Theta^C_\pi:C^\p(F)\to \BC$ be the induced function on the space $C^\p(F)$ of strongly regular semisimple conjugacy classes of $G(F)$. The trace of $h\in \CH$ can now be written as an integral over $C(F)$
\begin{equation}\label{eq:CC-OI:1}
\tr(\pi(h)) = \int_{C(F)}
\RO(h)(a)
\Theta^C_\pi(a)\ome_C = \langle \RO(h),\Theta^C_\pi\rangle
\end{equation}
where $\ome_C$ is the measure on $C(F)$ introduced in the last subsection, and $\langle \cdot,\cdot\rangle$ is the standard pairing of functions on $C(F)$ with respect to the measure $\ome_C$.

Similarly, for any local $L$-parameter $\vphi$ of $G(F)$ with associated $L$-packet $\Pi_\vphi$, let $\Theta_{\vphi}$ be the stable distribution character attached to $\vphi$ \cite{MR0540901}. Then $\Theta_\vphi$ descends to a function $\Theta^\Fc_\vphi$ on $\Fc(F)$ and the following identity holds
\begin{equation}\label{eq:CC-OI:2}
\int_{G(F)}
h(g)\Theta_\vphi(g)\ud g = 
\int_{\Fc(F)}
\RS\RO(h)(c)
\Theta^\Fc_{\vphi}(c)\ome_\Fc = 
\langle \RS\RO(h),\Theta^\Fc_{\vphi}\rangle.
\end{equation}

Finally, we discuss the construction of the Kazhdan cocenter, which can be viewed as a completion of the usual cocenter.

We first assume that $F$ is non-archimedean. By construction, the Bernstein center $\CZ$ acts on the Hecke algebra $\CH$ preserving the commutator $[\CH,\CH]$. Hence it acts on the cocenter $\RO(\CH) = \CH/[\CH,\CH]$. For any finite subset $I\subset \pi_0(\Ome)$, let $\CH_I = \CH \eps_I = \eps_I\CH$, which is a $\CZ_I =\CZ\eps_I$-module. The commutator $[\CH_I,\CH_I]$ is a $\CZ_I$-submodule, and hence $\CH_I/[\CH_I,\CH_I]$ is a $\CZ_I$-module. If $I^\p$ is a finite subset of $\pi_0(\Ome)$ containing $I$, then $\eps_{I^\p}\eps_I = \eps_I$. Thus we have the inclusion $\CH_I\to \CH_{I^\p}$ as well as a surjective $\CZ$-linear map $\CH_{I^\p}\to \CH_I$ given by $f\mapsto \eps_I f$. With respect to the inclusion, we have $[\CH_{I^\p},\CH_{I^\p}]\cap \CH_I = [\CH_I,\CH_I]$ \cite[Lem.~3.2]{KazCus}, and the surjection $f\mapsto \eps_I f$ sends $[\CH_{I^\p},\CH_{I^\p}]$ onto $[\CH_I,\CH_I]$. We derive an injective map $\CH_I/[\CH_I,\CH_I]\to \CH_{I^\p}/[\CH_{I^\p},\CH_{I^\p}]$ from the inclusion, and a surjective map $\CH_{I^\p}/[\CH_{I^\p},\CH_{I^\p}]\to \CH_I/[\CH_I,\CH_I]$ from the surjection. The inductive limit of $\CH_I/[\CH_I,\CH_I]$ over finite subsets $I$ is again the usual cocenter $\CH/[\CH,\CH]$. The corresponding projective limit is the \emph{Kazhdan cocenter}, denoted $\CC=\CC_G$, which is a completion of the usual cocenter $\CH/[\CH,\CH]$. Indeed, we have
$$
\CH/[\CH,\CH] = \CC^\fin
$$
where $\CC^\fin$ is the $\CZ$-submodule consisting of elements $f\in \CC$ such that $\eps_I f = f$ for some finite set $I$ of $\pi_0(\Ome)$.

The trace Paley-Wiener theorem \cite{MR874050} gives a spectral description of $\CC^\fin$. A linear form $\ell$ on the vector space $R(G(F))$ with basis $\Irr(G(F))$ is of the form
$$
\ell(\pi) = \langle \RO(h),\Theta^C_\pi\rangle
$$
for some $h\in \CH$ if and only if it satisfies the following two conditions: for every Levi subgroup $M$ of $G$ and every irreducible admissible representation $\sig$ of $M(F)$, the function $\chi\mapsto\ell(\Ind_P^G(\sig\otimes\chi))$ is regular on the torus of unramified characters of $M(F)$, where $P$ is a parabolic subgroup with Levi component $M$ and induction is normalized; and there exists an open compact subgroup $K\subset G(F)$ such that $\ell(\pi)=0$ whenever $\pi^K=0$. Such linear forms are called good forms. The Kazhdan cocenter $\CC$ is obtained by taking compatible good forms on finite unions of Bernstein components, without a global finiteness requirement. This characterization should be distinguished from unrestricted regular functions on the Lafforgue variety of \cite{psaromiligkos_lafforgue_2023}: at reducible induction parameters, trace values must satisfy the additive relations among the irreducible constituents.

\subsection{Archimedean variant}

Let $G$ be a reductive group over $F=\BR$. One can define the archimedean Kazhdan cocenter to be the projective limit over finite subsets of the set of $K$-types of $G(F)$. Precisely, let $\hat{K}$ be the set of equivalence classes of irreducible finite-dimensional representations of a fixed maximal compact subgroup of $G(F)$. For any $\sig\in \hat{K}$, let $e_\sig$ be the idempotent attached to $\sig$. For any finite subset $I$ of $\hat{K}$, let $\eps_I = \sum_{\sig\in I}e_\sig$. Following the same idea as in the non-archimedean case, we define the Kazhdan cocenter as the projective limit, over finite subsets $I$ of $\hat{K}$, of the level $I$ cocenters $\CS_I/[\CS_I,\CS_I]$, where $\CS_I = \eps_I\CS \eps_I$ and both the quotient $\CS_I/[\CS_I,\CS_I]$ and the projective limit are taken in the category of Fréchet spaces.

Alternatively, one can use the trace Paley-Wiener theorem for the $K$-finite Schwartz algebra $\CS_{(K)}$ and then drop the $K$-finiteness condition in its spectral characterization. This gives another definition, or characterization, of the archimedean analogue of the Kazhdan cocenter. Following \cite[Théorème~1.9]{MR0860667}, suppose $P=MAN$ is a cuspidal parabolic subgroup of $G$ with Langlands decomposition. Let $\hat{M}_{\mathrm{l.s.d}}$ be the set of irreducible representations of $M(F)$ that are non-degenerate limits of discrete series \cite[\S 12]{MR672840}. For $\del\in \hat{M}_{\mathrm{l.s.d}}$ and $\lam\in \Fa^*_\BC$, where $\Fa = \Lie A$, let $\pi_{\del,\lam}$ be the normalized representation induced from $P=MAN$ to $G$ from $\del\otimes e^\lam\otimes 1_N$.

Let $P_j = M_jA_jN_j$ ($j=1,\ldots,S$) be the collection of cuspidal parabolic subgroups of $G$ containing the fixed minimal parabolic of $G$. Let $\RW(A_j)$ be the Weyl group of $A_j$ that acts on $(\hat{M}_j)_{\mathrm{l.s.d}}$ and $\Fa^*_{j,\BC}$. The following theorem is established in \cite[Théorème~1.9]{MR0860667}.
\begin{thm}\label{thm:delorme-K-finite-Schwartz}
Let
$$
F_j:
\big(
\hat{M}_j
\big)_{\mathrm{l.s.d}}
\times \Fa^*_{j,\BC}\to \BC,\quad j=1,\ldots,S.
$$
The following two conditions are equivalent:
\begin{enumerate}
\item There exists $\phi\in \CS_{(K)}$ such that $F_j(\del,\lam) = \tr\pi_{\del,\lam}(\phi)$ for any $(\del,\lam)\in (\hat{M}_j)_{\mathrm{l.s.d}}\times \Fa^*_{j,\BC}$, $j=1,\ldots,S$.

\item The functions $F_j$ enjoy the following properties:
\begin{enumerate}
\item For any $j$, $F_j$ is finitely supported in $\del$.

\item For any $\del$, $F_j(\del,\cdot)\in \CF\CC^0(\Fa_j)$, where $\CF\CC^0(\Fa_j)$ consists of entire functions on $\Fa_{j,\BC}$ that are of moderate growth in bounded vertical strips. For details, see \cite[p.696]{MR0860667};

\item For any $j$, $F_j(w\del,w\lam) = F_j(\del,\lam)$ for any $w\in \RW(A_j)$;

\item Suppose that $P_i\subset P_j$. Then $M_i\subset M_j$ and $A_i\supset A_j$. Denote by $P_i\cap M_j = M_iA^j_iN^j_i$ the Langlands decomposition of the parabolic $P_i\cap M_j$ of $M_j$. Then $A^j_i = A_i\cap M_j$ and $N^j_i = N_i\cap M_j$. The decomposition $\Fa_i = \Fa_j\oplus \Fa^j_i$ identifies $\Fa^*_j$ with a subspace of $\Fa^*_i$. Suppose that $\Ind^{M_j}_{P_i\cap M_j}(\del_i\otimes e^0\otimes 1) = \bigoplus_{l=1}^R \del_{j,l}$ with $\del_i\in (\hat{M}_i)_{\mathrm{l.s.d}}$ and $\del_{j,l}\in (\hat{M}_j)_{\mathrm{l.s.d}}$. Then
$$
\forall \lam\in \Fa^*_{j,\BC}\subset \Fa^*_{i,\BC}, \quad 
F_i(\del_i,\lam) = 
\sum_{l=1}^R F_j(\del_{j,l},\lam).
$$
\end{enumerate}
\end{enumerate}
\end{thm}
Using the preceding theorem, we may drop the finiteness assumption in Part (2)(a) of Theorem \ref{thm:delorme-K-finite-Schwartz} and introduce the following spectral characterization of the Kazhdan cocenter in the archimedean setting.

\begin{defin}\label{defin:KCC-archimedean}
Following the notation in Theorem \ref{thm:delorme-K-finite-Schwartz}, the Kazhdan cocenter $\CC$ is defined to be the collection of functions 
$$
F_j:
\big(
\hat{M}_j
\big)_{\mathrm{l.s.d}}
\times \Fa^*_{j,\BC}\to \BC,\quad j=1,\ldots,S
$$
with the following properties:
\begin{enumerate}
\item For any $\del$, $F_j(\del,\cdot)\in \CF\CC^0(\Fa_j)$, where $\CF\CC^0(\Fa_j)$ consists of entire functions on $\Fa_{j,\BC}$ that are of moderate growth in bounded vertical strips. For details, see \cite[p.696]{MR0860667};

\item For any $j$, $F_j(w\del,w\lam) = F_j(\del,\lam)$ for any $w\in \RW(A_j)$;

\item Suppose that $P_i\subset P_j$. Then $M_i\subset M_j$ and $A_i\supset A_j$. Denote by $P_i\cap M_j = M_iA^j_iN^j_i$ the Langlands decomposition of the parabolic $P_i\cap M_j$ of $M_j$. Then $A^j_i = A_i\cap M_j$ and $N^j_i = N_i\cap M_j$. The decomposition $\Fa_i = \Fa_j\oplus \Fa^j_i$ identifies $\Fa^*_j$ with a subspace of $\Fa^*_i$. Suppose that $\Ind^{M_j}_{P_i\cap M_j}(\del_i\otimes e^0\otimes 1) = \bigoplus_{l=1}^R \del_{j,l}$ with $\del_i\in (\hat{M}_i)_{\mathrm{l.s.d}}$ and $\del_{j,l}\in (\hat{M}_j)_{\mathrm{l.s.d}}$. Then
$$
\forall \lam\in \Fa^*_{j,\BC}\subset \Fa^*_{i,\BC}, \quad 
F_i(\del_i,\lam) = 
\sum_{l=1}^R F_j(\del_{j,l},\lam).
$$
\end{enumerate}
\end{defin}

As with the Bernstein center, we can introduce the finite (resp. stable) variant of the
Kazhdan cocenter $\CC^{\fin}$ (resp. $\CC^\st$) and $\CC^{\st,\fin}$.
We leave the details to the reader.

\section{Lafforgue transform}\label{sec:lafforgue_transform}

In this section, we construct the Lafforgue transform, a natural map from the
Kazhdan cocenter to the Bernstein center. The construction is representation-theoretic
in nature: it converts cocentral, or orbital, data into spectral multipliers. The material
reviewed here is mainly taken from \cite[\S 5]{LN}. We state the main results and explain
their meaning, referring the proofs to \cite{LN}.

\subsection{Definition and first properties}\label{subsec:Laf-definition}

Fix a reductive group $G$ over a non-archimedean local field $F$. Let
\[
\Irr=\Irr_G=\Irr(G(F))
\]
be the set of irreducible smooth admissible representations of $G(F)$. By the work of
Psaromiligkos \cite{psaromiligkos_lafforgue_2023}, this set is parametrized, block by block,
by the open simple-module locus in the corresponding \emph{Lafforgue variety}. The
cuspidal-support map
\[
\Irr(G(F))\longrightarrow \Ome(G(F))
\]
is algebraic and finite-to-one. The Bernstein center $\CZ=\CZ_G$ is
identified with the ring of regular functions on $\Ome(G(F))$.

On the other hand, every element of the Kazhdan cocenter $\CC=\CC_G$ gives a trace
function on $\Irr(G(F))$, with the regularity in induced families described above. Indeed, for $c\in \CC$ and
$\pi\in \Irr(G(F))$, one can evaluate the trace
\[
\Tr(c,\pi).
\]
When $c$ is represented by a test function, this is the usual trace of the operator on
$\pi$, and the trace only depends on the image of the test function in the cocenter.

There exists an open dense subvariety
\[
\Ome^\p(G)\subset \Ome(G)
\]
consisting of cuspidal data whose parabolic induction to $G(F)$ is irreducible. Its inverse
image
\[
\Irr^\p(G(F))\subset \Irr(G(F))
\]
is isomorphic to $\Ome^\p(G)$. For every $c\in \CC$, its trace function on
$\Irr^\p(G(F))\simeq \Ome^\p(G)$ extends uniquely to a regular function on the whole
Bernstein variety: the extension is given by traces on the full parabolically induced
representations, whose regularity follows from the trace Paley-Wiener theorem.
This gives the Lafforgue transform
\begin{equation}\label{eq:Laf-Rep:1}
\Laf:\CC\longrightarrow \CZ.
\end{equation}

Equivalently, if $[M,\sig]$ is a cuspidal datum, where $M$ is a Levi subgroup of $G$ and
$\sig$ is a supercuspidal representation of $M(F)$, then
\begin{equation}\label{eq:Laf-Rep:2}
\Laf(c)\big([M,\sig]\big)=
\Tr\big(c,\Ind(\sig)\big),
\end{equation}
where $\Ind$ denotes normalized parabolic induction from $M(F)$ to $G(F)$.

\begin{pro}[{\cite[\S 5.1]{LN}}]\label{pro:Laf-surjective}
The Lafforgue transform
\[
\Laf:\CC\to \CZ
\]
is a $\CZ$-linear surjective map. After restricting to the finite part, it induces a
$\CZ$-linear morphism
\[
\Laf:\CC^{\fin}\to \CZ^{\fin}.
\]
Moreover, $\Laf$ sends $\CC^{\st,\fin}$ into $\CZ^{\st,\fin}$.
\end{pro}

The last assertion should be interpreted as follows. Elements in $\CC^{\st,\fin}$ are
represented by stable orbital integrals of test functions. On the other hand, distributions
in $\CZ^\fin$ are represented by invariant locally integrable functions on $G(F)$
\cite[Thm.~2.5]{MR1979109}. It follows that the elements in
$\Laf(\CC^{\st,\fin})$ are represented by invariant locally integrable functions on $G(F)$
that are also stably invariant, and hence by stably invariant functions on $G(F)$.

\subsection{Compatibility with stable transfer}\label{subsec:Laf-stable-transfer}

The Lafforgue transform is expected to be compatible with Langlands stable transfer. Suppose
that
\[
\rho:{}^LH\to {}^LG
\]
is a homomorphism of $L$-groups. Langlands stable transfer should give a map
\[
\CT^\st_\rho:\CC_G^{\st}\to \CC_H^{\st}
\]
dual to the transfer of stable characters. On the spectral side, one expects a homomorphism
\[
\CZ^{\st}_\rho:\CZ_G^{\st}\to \CZ_H^{\st}
\]
induced by the corresponding functorial map of $L$-packets. The expected compatibility is
expressed by the commutativity of the following diagram:
\begin{equation}\label{eq:Laf-stable-transfer-diagram}
\xymatrix{
\CC_G^{\st} \ar[r]^{\CT^\st_\rho} \ar[d]_{\Laf_G^{\st}} &
\CC_H^{\st} \ar[d]^{\Laf_H^{\st}} \\
\CZ_G^{\st} \ar[r]_{\CZ^{\st}_\rho} &
\CZ_H^{\st}.
}
\end{equation}

This expectation explains why the Lafforgue transform is natural in the
Braverman-Kazhdan setting. Stable transfer is naturally formulated on the stable cocenter,
whereas local gamma factors live on the Bernstein center side. Thus the Lafforgue transform
is the mechanism that connects the orbital side with the spectral side.

\subsection{Sections of the Lafforgue transform}\label{subsec:sections-Laf}

In general, $\Laf$ is not injective, and therefore there are many possible choices of sections.

For the applications in \cite{LN}, one uses a distinguished section with the following properties. Let
$G=\SL_2$ or $\GL_2$. Let $\CC^{\perp,\triv}$ be the subspace of $\CC$ consisting of
elements orthogonal to one-dimensional characters of $G(F)$, namely
\[
\langle c,\chi\rangle=0
\]
for every one-dimensional character $\chi$ of $G(F)$. There exists a section
\[
\LafSec:\CZ\to \CC^{\perp,\triv}\subset \CC
\]
such that
\[
\LafSec(z)\big(\Ind_M^G(\sig)\big)=z\big([M,\sig]\big)
\]
for every cuspidal datum $[M,\sig]$.

When restricted to the finite stable parts, elements in $\CZ^{\st,\fin}$ and
$\CC^{\st,\fin}$ can be realized as functions on the Steinberg-Hitchin base $\Fc(F)$.
One of the main results of \cite{LN} is that, for $G=\SL_2$ or $\GL_2$, the section
$\LafSec$ and the transform $\Laf$ can be realized as explicit integral transforms on
$\Fc(F)$.

\begin{thm}[{\cite[Thm.~5.1.3]{LN}}]\label{thm:stable-Laf}
Let $G=\SL_2$ or $\GL_2$ over a local field $F$ of residual characteristic not equal to two.
For any distribution $z\in \CZ^{\st,\fin}$ lying in the stable Bernstein center and supported
on finitely many Bernstein components, there exists a stably invariant function
\[
\CJ_G=\CJ_\Fc\circ \bc
\]
such that $z=\CJ_G\ud g$. Then $\LafSec(z)$ is given, as a function on the
Steinberg-Hitchin base, by the formula
\begin{equation}\label{eq:Laf-formula}
\LafSec(z)=
\CF_\psi\left(
\frac{\CF_{\psi^{-1}}(\CJ_\Fc)}{|\cdot|}
\right),
\end{equation}
where the Fourier transforms and $|\cdot|$ are taken in the trace variable.
\end{thm}

Theorem \ref{thm:stable-Laf} remains valid when $F$ has residual characteristic two; this extension will be presented in the forthcoming joint work \cite{NL-general}.

In the range treated in \cite{LN}, the proof of Theorem \ref{thm:stable-Laf} relies on
the Gelfand-Graev formula for stable transfer and the descent formula for dihedral
representations; see \cite[\S 5--6]{LN}. The
important point for the present report is that the representation-theoretic section
$\LafSec$ is realized by an explicit Fourier-type integral transform on the Steinberg-Hitchin
base, with a formula of motivic nature valid over every local field.

The inverse transform on the stable finite part is formally given by
\begin{equation}\label{eq:Laf-inverse-formula}
\Laf(c)=
\CF_\psi\big(
|\cdot|\,\CF_{\psi^{-1}}(c)
\big),
\end{equation}
whenever $c$ lies in the corresponding stable finite cocenter and the expression is defined.
The analytic meaning of these formulas is justified in \cite[\S 5.5]{LN} using explicit
calculations of Moy and Tadic.

\subsection{Supercuspidal idempotents}\label{subsec:Laf-supercuspidal}

We recall the effect of $\LafSec$ on supercuspidal idempotents for $G=\SL_2$, following
\cite[\S 5.2]{LN}. Let $\pi$ be a supercuspidal representation of $G(F)$. It gives an
isolated point of the Bernstein variety, and we denote by $e_\pi\in \CZ$ the idempotent
function on $\Ome$ that takes the value one at $\pi$ and zero on the other connected
components.

\begin{thm}[{\cite[Thm.~5.2.1]{LN}}]\label{thm:Laf-supercuspidal}
Let $G=\SL_2$ and let $\pi$ be a supercuspidal representation of $G(F)$.
\begin{enumerate}
\item The idempotent $e_\pi$ is represented by
\[
e_\pi=\Theta_\pi d_\pi\ud g,
\]
where $\Theta_\pi$ is the distribution character of $\pi$ and $d_\pi$ is the formal degree.
The formal degree depends on the choice of Haar measure, but the product $d_\pi\ud g$ is
independent of that choice.

\item As a function on the set of regular semisimple conjugacy classes, one has
\[
\LafSec(e_\pi)=\Delta\,v\,\mathbf{1}_{\el}\,\Theta_\pi.
\]
Here $\Delta$ is the square root of the absolute value of the Weyl discriminant,
$\mathbf{1}_{\el}$ is the characteristic function of the elliptic regular semisimple locus,
and $v$ is the inverse of the volume function on the elliptic regular semisimple locus.
\end{enumerate}
\end{thm}

Thus the Lafforgue section sends the supercuspidal Bernstein idempotent, represented
spectrally by the character times the formal degree, to the normalized elliptic character
in the cocenter. This is one of the clearest concrete manifestations of the philosophy that
the Lafforgue transform connects the spectral and orbital sides.

\subsection{Real Chebyshev formula}\label{subsec:Laf-Chebyshev}

For discrete series representations of $\SL_2(\BR)$, there is an archimedean analogue of
Theorem \ref{thm:Laf-supercuspidal}. Let $\Theta_\ell$ be the stable discrete series
character of weight $\ell\geq 1$, and let $d_\ell=\ell$ be its formal degree. The following
result is proved in \cite[\S 5.3]{LN}.

\begin{thm}[{\cite[Thm.~5.3.1]{LN}}]\label{thm:Laf-Chebyshev}
For $G=\SL_2(\BR)$, the Lafforgue transform sends $d_\ell\Theta_\ell$ to the elliptic part
$v\Delta\Theta_\ell^{\el}$. Equivalently,
\begin{equation}\label{eq:Laf-Chebyshev}
\CF_\psi\left(
\frac{\CF_{\psi^{-1}}(d_\ell\Theta_\ell)}{|\cdot|}
\right)
=
v\Delta\Theta_\ell^{\el}.
\end{equation}
\end{thm}

The proof is an elementary calculation with Chebyshev polynomials. Writing
\[
x=e^{i\theta}+e^{-i\theta}=2\cos\theta,
\]
one has
\[
\frac{e^{i\ell\theta}+e^{-i\ell\theta}}{|e^{i\theta}-e^{-i\theta}|}
=
\frac{2T_\ell(x/2)}{\sqrt{4-x^2}},
\]
where $T_\ell$ is the Chebyshev polynomial of the first kind. The elliptic discrete series
character is expressed in terms of the Chebyshev polynomial of the second kind. The identity
in Theorem \ref{thm:Laf-Chebyshev} is then reduced to the standard differential relation
between the two families of Chebyshev polynomials.

\subsection{The image of the stable Lafforgue transform}\label{subsec:image-stable-Laf}

For $G=\SL_2$ over a non-archimedean local field of odd residual characteristic, the image
of the stable Lafforgue section admits a simple description.

\begin{thm}[{\cite[Thm.~5.4.1]{LN}}]\label{thm:image-stable-Laf}
Let $G=\SL_2$ over $F$, where $F$ has odd residual characteristic. The image of
\[
\LafSec:\CZ^{\st}\to \CC^{\st}
\]
is the subspace $\CC^{\st,\perp,\triv}$ consisting of elements of the stable Kazhdan cocenter
that are orthogonal to the character of the trivial representation.
\end{thm}

We recall the idea of the proof from \cite[\S 5.4]{LN}. First, the analytic vanishing
property in Theorem \ref{thm:analytic-J} below implies that the image of $\LafSec$ is
orthogonal to the trivial character. Second, the trace Paley-Wiener theorem and the
description of stable packets for $\SL_2$ show that every element of
$\CC^{\st,\perp,\triv}$ is obtained from the stable Bernstein center by $\LafSec$. Similar
statements hold for $G=\GL_2$, with the expected modifications involving central characters.

\subsection{Analytic properties of the stable Bernstein-center functions}
\label{subsec:analytic-J}

The formula \eqref{eq:Laf-formula} is not purely formal. One must know that the quotient
by $|\cdot|$ after Fourier transform is integrable. This is supplied by the following theorem
of \cite[\S 5.5]{LN}, whose proof uses explicit formulas of Moy and Tadic in the
non-archimedean case.

\begin{thm}[{\cite[Thm.~5.5.1]{LN}}]\label{thm:analytic-J}
Let $G=\SL_2$ or $\GL_2$, and assume that the residual characteristic of $F$ is not equal
to two. For any $\CJ_G\in \CZ_G^{\st,\fin}$, the following statements hold.
\begin{enumerate}
\item There exists an absolutely integrable function $\CJ_\Fc$ on the Steinberg-Hitchin base
$\Fc_G(F)$ such that
\[
\CJ_G=(\CJ_\Fc\circ \bc)\ud g.
\]
For $G=\GL_2$, the function $\CJ$ is rapidly decreasing in the determinant variable.

\item The inverse Fourier transform $\CF_{\psi^{-1}}(\CJ_\Fc)$ in the trace variable is represented
by a globally bounded smooth function that vanishes at zero.

\item The function
\[
\frac{\CF_{\psi^{-1}}(\CJ_\Fc)}{|\cdot|}
\]
is absolutely integrable in the trace variable. Consequently, for any \'{e}tale quadratic
$F$-algebra $K$, the function
\[
\CF_\psi\left(
\frac{\eta_K\cdot \CF_{\psi^{-1}}(\CJ_\Fc)}{|\cdot|}
\right)
\]
is globally bounded, where $\eta_K$ is the quadratic character of $F^\times$ determined by $K$.
\end{enumerate}
\end{thm}

This theorem has two consequences for the present discussion. First, it justifies the integral
formula \eqref{eq:Laf-formula} defining $\LafSec(z)$ on the Steinberg-Hitchin base. Second,
the vanishing of $\CF_{\psi^{-1}}(\CJ_\Fc)$ at zero explains why the image of the stable
Lafforgue section is orthogonal to the trivial character, as stated in Theorem
\ref{thm:image-stable-Laf}.

\begin{rmk}
Since orbital integrals normalized by the square root of the Weyl discriminant are locally
bounded, the first assertion of Theorem \ref{thm:analytic-J} implies that the corresponding
distribution $\RJ_G$ is locally integrable on $G(F)$. For $G=\GL_2$, the rapid decay in the
determinant variable is obtained by applying Mellin transform along the center.
\end{rmk}

\subsection{Further expectations}\label{subsec:Laf-expectations}

For a general reductive group $G$ over a local field $F$, one expects a
distinguished section
\[
\LafSec:\CZ\to \CC
\]
of the Lafforgue transform. After restriction to $\CZ^{\st,\fin}$, both $\Laf$ and
$\LafSec$ should be expressible by explicit integral transforms on the Steinberg-Hitchin
base of $G$. In rank one, Theorem \ref{thm:stable-Laf} gives such a formula. For general
reductive groups, Ng\^o and the author have formulated conjectural analytic formulas
for the stable Lafforgue transform and its section on the Steinberg-Hitchin base.
These formulas will be presented in the work in preparation \cite{NL-general}.

One also expects the image of the distinguished section to satisfy a natural orthogonality
condition. In the rank-one cases considered above, the image is orthogonal to one-dimensional
characters. More generally, the expected image should be orthogonal to the distribution
characters of irreducible smooth admissible representations of $G(F)$ that are not essentially
tempered.

\newcommand{\etalchar}[1]{$^{#1}$}
\providecommand{\bysame}{\leavevmode\hbox to3em{\hrulefill}\thinspace}
\providecommand{\MR}{\relax\ifhmode\unskip\space\fi MR }
\providecommand{\MRhref}[2]{%
  \href{http://www.ams.org/mathscinet-getitem?mr=#1}{#2}
}
\providecommand{\href}[2]{#2}


\begin{thebibliography}{GPSR87b}

\bibitem[BDK86]{MR874050}
J.~Bernstein, P.~Deligne, and D.~Kazhdan, \emph{Trace {P}aley-{W}iener theorem
  for reductive {$p$}-adic groups}, J. Analyse Math. \textbf{47} (1986),
  180--192.

\bibitem[Ber84a]{BD84}
J.~Bernstein, \emph{Le ``centre'' de {B}ernstein}, Representations of reductive
  groups over a local field, Travaux en Cours, Hermann, Paris, 1984, Edited by
  P. Deligne, pp.~1--32.

\bibitem[BK00]{BK00}
A.~Braverman and D.~Kazhdan, \emph{{$\gamma$}-functions of representations and
  lifting}, Geom. Funct. Anal. (2000), no.~Special Volume, Part I, 237--278,
  With an appendix by V. Vologodsky, GAFA 2000 (Tel Aviv, 1999).

\bibitem[BK14]{bkglobalization}
J.~Bernstein and B.~Kr\"{o}tz, \emph{Smooth {F}r\'{e}chet globalizations of
  {H}arish-{C}handra modules}, Israel J. Math. \textbf{199} (2014), no.~1,
  45--111.

\bibitem[BKV13]{bezrukavnikov2013categorical}
R.~Bezrukavnikov, D.~Kazhdan, and Y.~Varshavsky, \emph{A categorical
  approach to the stable center conjecture}, arXiv preprint arXiv:1307.4669
  (2013).

\bibitem[Del86]{MR0860667}
P. Delorme, \emph{Formules limites et formules asymptotiques pour les
  multiplicit\'{e}s dans {$L^2(G/\Gamma)$}}, Duke Math. J. \textbf{53} (1986),
  no.~3, 691--731.

\bibitem[Del05]{MR2183287}
\bysame, \emph{Sur le th\'{e}or\`eme de {P}aley-{W}iener d'{A}rthur}, Ann. of
  Math. (2) \textbf{162} (2005), no.~2, 987--1029.

\bibitem[FLN10a]{FLN10}
E.~Frenkel, R.~Langlands, and B.~C. Ng{\^o}, \emph{Formule des traces et
  fonctorialit\'{e}: le d\'{e}but d'un programme}, Ann. Sci. Math. Qu\'{e}bec
  \textbf{34} (2010), no.~2, 199--243.

\bibitem[GG63]{MR0155931}
I.~M. Gel'fand and M.~I. Graev, \emph{Representations of the group of
  second-order matrices with elements in a locally compact field and special
  functions on locally compact fields}, Uspehi Mat. Nauk \textbf{18} (1963),
  no.~4(112), 29--99.

\bibitem[GJ72]{gjzeta}
R.~Godement and H.~Jacquet, \emph{Zeta functions of simple algebras}, Lecture
  Notes in Mathematics, Vol. 260, Springer-Verlag, Berlin-New York, 1972.

\bibitem[Hai14]{MR3444233}
T.~Haines, \emph{The stable {B}ernstein center and test functions for
  {S}himura varieties}, Automorphic forms and {G}alois representations. {V}ol.
  2, London Math. Soc. Lecture Note Ser., vol. 415, Cambridge Univ. Press,
  Cambridge, 2014, pp.~118--186.

\bibitem[HC70]{HC70}
Harish-Chandra, \emph{Harmonic analysis on reductive {$p$}-adic groups}, Lecture Notes
  in Mathematics, Vol. 162, Springer-Verlag, Berlin-New York, 1970, Notes by G.
  van Dijk.

\bibitem[JL21]{jiang2021certain}
D.~Jiang and Z.~Luo, \emph{Certain {F}ourier operators and their associated
  poisson summation formulae on $\mathrm{GL}_1$}, arXiv:2108.03566 (2021).

\bibitem[JL22]{MR4474366}
\bysame, \emph{Certain {F}ourier operators on {${\rm GL}_1$} and local
  {L}anglands gamma functions}, Pacific J. Math. \textbf{318} (2022), no.~2,
  339--374.

\bibitem[JL23]{DanielZhilin}
D.~Johnstone and Z.~Luo, \emph{Stable transfer for functorial lifts of
  $\mathrm{GL}_2$ over local fields}, preprint (2023).

\bibitem[JLZ22]{JLZ}
D.~Jiang, Z.~Luo, and L.~Zhang, \emph{Harmonic analysis and gamma functions on
  symplectic groups}, To appear in Memoirs of the American Mathematical Society
  (2022).

\bibitem[Kaz86a]{KazCus}
D.~Kazhdan, \emph{Cuspidal geometry of {$p$}-adic groups}, J. Analyse Math.
  \textbf{47} (1986), 1--36.

\bibitem[Kot82]{kott2}
R.~Kottwitz, \emph{Rational conjugacy classes in reductive groups}, Duke
  Math. J. \textbf{49} (1982), no.~4, 785--806.

\bibitem[Laf14a]{lafforguegl2}
L.~Lafforgue, \emph{Expos{\'e} 3 (jeudi 3 juillet) : Le cas de $\mathrm{GL}(2)$
  : {\'e}tude locale encore formelle (avec un paragraphe suppl{\'e}mentaire sur
  le cas g{\'e}n{\'e}ral)}.

\bibitem[Laf16]{lafforgue2016principe}
\bysame, \emph{Le principe de fonctorialit{\'e} de {L}anglands comme un
  probl{\`e}me de g{\'e}n{\'e}ralisation de la loi d'addition}, Tech. report,
  IHES preprint,, 2016.

\bibitem[Lan70]{langlandsproblems}
R.~Langlands, \emph{Problems in the theory of automorphic forms}, 18--61.
  Lecture Notes in Math., Vol. 170.

\bibitem[Lan79]{MR0540901}
\bysame, \emph{Stable conjugacy: definitions and lemmas}, Canadian J.
  Math. \textbf{31} (1979), no.~4, 700--725.

\bibitem[Lan13]{MR3117742}
\bysame, \emph{Singularit\'{e}s et transfert}, Ann. Math. Qu\'{e}.
  \textbf{37} (2013), no.~2, 173--253.

\bibitem[LN]{LN}
Z.~Luo and B.~C. Ng\^{o}, \emph{Nonabelian {F}ourier kernels on $\mathrm{SL}_2$
  and $\mathrm{GL}_2$}, arXiv:2409.14696.

\bibitem[LN2]{NL-general}
Z.~Luo and B.~C. Ng\^{o}, \emph{Lafforgue transform for general reductive groups},
  in preparation.

\bibitem[Luo]{Zhilinsingapore}
Z.~Luo, \emph{An introduction to the proposal of {B}raverman and {K}azhdan},
  Singapore Lecture Notes.

\bibitem[Luo19a]{MR3990815}
\bysame, \emph{On the {B}raverman-{K}azhdan proposal for local factors:
  spherical case}, Pacific J. Math. \textbf{300} (2019), no.~2, 431--471.

\bibitem[MT02]{MR1979109}
A. Moy and M. Tadi\'{c}, \emph{The {B}ernstein center in terms of
  invariant locally integrable functions}, Represent. Theory \textbf{6} (2002),
  313--329.

\bibitem[Ng{\^o}20]{ngo2016hankel}
B.~C. Ng{\^o}, \emph{Hankel transform, {L}anglands functoriality and functional
  equation of automorphic {$L$}-functions}, Jpn. J. Math. \textbf{15} (2020),
  no.~1, 121--167.

\bibitem[SS13]{MR2983012}
P. Scholze and S. Shin, \emph{On the cohomology of compact unitary
  group {S}himura varieties at ramified split places}, J. Amer. Math. Soc.
  \textbf{26} (2013), no.~1, 261--294.

\bibitem[Ste74]{Steinberg-conjugacy}
R.~Steinberg, \emph{Conjugacy classes in algebraic groups}, Lecture Notes in
  Mathematics, vol. 366, Springer-Verlag, Berlin-New York, 1974, Notes by Vinay
  V. Deodhar.

\bibitem[Tat50]{tatethesis}
J.~Tate, \emph{Fourier analysis in number fields and {H}ecke's zeta-functions},
  ProQuest LLC, Ann Arbor, MI, 1950, Thesis (Ph.D.)--Princeton University.

\bibitem[Wal88]{Wallach-RealI}
N.~R. Wallach, \emph{Real reductive groups. {I}}, Pure and Applied Mathematics,
  vol. 132, Academic Press, Inc., Boston, MA, 1988.

\bibitem[BPLZZ21]{MR4298750}
R.~Beuzart-Plessis, Y.~Liu, W.~Zhang, and X.~Zhu, \emph{Isolation of cuspidal
  spectrum, with application to the {G}an-{G}ross-{P}rasad conjecture}, Ann. of
  Math. (2) \textbf{194} (2021), no.~2, 519--584.

\bibitem[BT65]{MR207712}
A.~Borel and J.~Tits, \emph{Groupes r\'eductifs}, Inst. Hautes \'Etudes Sci.
  Publ. Math. (1965), no.~27, 55--150.

\bibitem[KZ82]{MR672840}
A.~Knapp and G.~Zuckerman, \emph{Classification of irreducible tempered
  representations of semisimple groups}, Ann. of Math. (2) \textbf{116} (1982),
  no.~2, 389--455.

\bibitem[Psa23b]{psaromiligkos_lafforgue_2023}
K.~I. Psaromiligkos, \emph{The {Lafforgue} variety and irreducibility of
  induced representations}, May 2023, arXiv:2211.11834 [math].

\bibitem[Var77]{varrealHC}
V.~S. Varadarajan, \emph{Harmonic analysis on real reductive groups}, Lecture
  Notes in Mathematics, vol. Vol. 576, Springer-Verlag, Berlin-New York, 1977.

\end{thebibliography}
\end{document}